\documentclass[mn,fleqn]{w-art}
\usepackage{times}
\usepackage{amsmath,amssymb}
\usepackage{w-thm}
\usepackage[]{graphicx}

\begin{document}
\pagespan{1}{}
\title[Osculating properties of decomposable scrolls]{Osculating properties of decomposable scrolls}
\author[Antonio Lanteri]{Antonio Lanteri\footnote{Corresponding author:
    e-mail: {\sf lanteri@mat.unimi.it}, Phone: +39\,02\,50316153, Fax:
    +39\,02\,50316090}\inst{1}} \address[\inst{1}]{Dipartimento di
Matematica ``F. Enriques'', Universit\`{a} degli Studi di
Milano, Via C. Saldini, 50, I-20133 Milano,
Italia}

\author[Raquel Mallavibarrena]{Raquel Mallavibarrena\footnote{
    e-mail: {\sf rmallavi@mat.ucm.es}, Phone: +34\,91\,3944657, Fax:
    +34\,91\,3944607}\inst{2}}
\address[\inst{2}]{Departamento de Algebra, Facultad
de Matem\'{a}ticas, Universidad Complutense de Madrid, Ciudad
Universitaria, E-28040 Madrid, Spain}



\begin{abstract} Osculating spaces of decomposable scrolls (of any genus
and not necessarily normal) are studied and their inflectional
loci are related to those of their generating curves by using
systematically an idea introduced by Piene and Sacchiero in the
setting of rational normal scrolls. In this broader setting the
extra components of the second discriminant locus --deriving from
flexes-- are investigated and a new class of uninflected surface
scrolls is presented and characterized. Further properties related
to osculation are discussed for (not necessarily decomposable)
scrolls.
\end{abstract}

\maketitle

\section*{Introduction}
The inflectional behavior of a projective variety belongs to its
extrinsic geometry. In particular, flexes can appear on projective
manifolds under (isomorphic) projections. Though this observation
is obvious, it seems that several projective manifolds have been
extensively investigated from the point of view of their
osculatory behavior only in the linearly normal case. This is true
e.g., for rational scrolls of any dimension \cite{bib9} and also
for elliptic surface scrolls \cite{bib7}. In this paper we mainly
consider {\it decomposable scrolls}, not necessarily linearly
normally embedded, and we study their inflectional behavior.
\par
Decomposable scrolls $X \subset \mathbb{P}^N$, whose construction
generalizes that of rational normal scrolls, are generated by $n$
curves $C_i$ ($i=1, \dots, n$) isomorphic each other, lying in
linearly independent linear subspaces generating the whole
$\mathbb{P}^N$ (see Section 1). They are very well suited to
investigate their $k$-th inflectional loci $\Phi_k(X)$. We do that
developing systematically the local description used in
\cite{bib9} and \cite{bib7}, and in Sections 1 and 2 we succeed to
describe several properties of $\Phi_k(X)$, relating them to the
inflectional loci of the generating sections $C_i$.
\par
In particular, restricting to the case of rational non-normal
scrolls our approach allows us to produce in Section 3 a new
series of counterexamples to the even dimensional part of a
conjecture of Piene and Tai \cite{bib10}. While the odd
dimensional part of this conjecture has been proved several years
ago \cite{bib10}, \cite{bib3}, the even dimensional part is false
for certain linearly normal scrolls, as shown by the first author
\cite{bib6}. However we want to stress that the new
counterexamples exhibited here are rational scrolls, though, of
course, not linearly normal. We also characterize these examples
in the framework of decomposable scrolls (Theorem 3.4). This adds
some information in order to correct the even dimensional part of
the conjecture.
\par
Let $X$ be a decomposable scroll. While describing $\Phi_k(X)$ for
$k > 2$ involves inflectional loci of lower order, the description
becomes very easy for $k=2$. In particular, we show that for a
decomposable scroll $X$, $\Phi_2(X)$ can have only two types of
irreducible components. Let $G$ be any such a component. Then,
either $G$ is a sub-fibre of a fibre of $X$, or $X$ is rational,
some curve $C_i$ is a line, and $G$ is a sub-scroll of $X$ given
by a Segre product (Proposition 4.2).
\par
This precise description of $\Phi_2(X)$ allows us to study in
Section 4 the second discriminant locus of a decomposable scroll
$X \subset \mathbb{P}^N$. This is the Zariski closed subset
$\mathcal{D}$ of $\mathbb{P}^{N \vee}$ parameterizing all
hyperplane sections of $X$ admitting a triple point. The main
component of $\mathcal{D}$ is the second dual variety of $X$,
which parameterizes osculating hyperplanes to $X$ at general
points and their limits. But when $X$ has flexes, extra components
$\mathcal{D}_G$ of $\mathcal{D}$ arise, coming from the
irreducible components $G$ of $\Phi_2(X)$. Our study of
$\Phi_2(X)$ allows us to describe these components: either
$\mathcal{D}_G$ is a linear space or it is a $1$-dimensional
family of linear spaces. In particular, we show that
$\mathcal{D}_G$ is a scroll if and only if $X$ is a rational
normal scroll generated by some lines plus conics and/or twisted
cubics (Example 4.3 and Proposition 4.4). Moreover, we
characterize rational normal scrolls generated by some lines plus
some conics as the decomposable scrolls admitting an irreducible
component $\mathcal{D}_G$ of $\mathcal{D}$ which is a rational
normal scroll (Theorem 4.7).
\par
In Section 5, we come to surface scrolls, not necessarily linearly
normal, regardless the fact they are decomposable or not. Here the
techniques developed in the previous sections fail. We discuss two
points arising from \cite{bib6}. a) Indecomposable elliptic
surface scrolls of invariant $-1$ have been studied in
\cite{bib6}. By adapting the approach used there, we investigate
those of invariant $0$, providing a description of their flexes in
terms of base points of suitable linear systems related to the one
giving the embedding (Proposition 5.1). b) The lowest dimension of
any osculating space to a surface scroll is $3$, as shown in
\cite{bib6}. Moreover, Example 3.2 shows that any higher order
osculating space can have very low dimension at some points. Here
we find sufficient conditions to grant that all $k$-th osculating
spaces of a surface scroll have dimension $\geq k+1$. They are
formulated in terms of the (relatively) good properties of the
linear system giving rise to the embedding (Theorem 5.2).

\setcounter{section}{-1}
\section{Notation and background}

We work over the field of complex numbers. Let $X$ be a smooth
projective variety of dimension $n \geq 1$. If $L$ is a line
bundle on $X$ we denote by $|W|$ the (not necessarily complete)
linear system defined by a vector subspace $W \subseteq H^0(X,L)$.
Suppose that $|W|$ is very ample, i.\ e., the map defined by $W$
is an embedding $\varphi_W:X \hookrightarrow
\mathbb{P}(W)=\mathbb{P}^N$. Then
$L=\varphi_W^*(\mathcal{O}_{\mathbb{P}^N}(1))$. In this case,
frequently we look at the pair $(X,W)$, or at the triplet
$(X,L,W)$, in place of the non-degenerate embedded variety
$\varphi_W(X) \subset \mathbb{P}^N$ and sometimes we do not
distinguish between $X$ and its image.
\par
For any integer $k \geq 0$ let $J_kL$ be the $k$-th jet bundle of
$L$. For every $x \in X$ we denote by
\[j_{k,x}^{(X,W)}:W \to (J_kL)_x\] the homomorphism associating to
every section $\sigma \in W$ its $k$-th jet evaluated at $x$. When
the subspace $W$ we are dealing with is clear from the context, or
the discussion involves a single pair $(X,W)$, we simply write
$j_{k,x}^X$ or $j_{k,x}$ respectively, instead of
$j_{k,x}^{(X,W)}$. Recall that $j_{k,x}(\sigma)$ is represented in
local coordinates by the Taylor expansion of $\sigma$ at $x$,
truncated after the order $k$. So, if $|W|$ is very ample, the
$k$-{\it{th osculating subspace}} to $X$ at a point $x \in X$ is
defined as
$\text{Osc}_x^k(X):=\mathbb{P}(\text{Im}j_{k,x}^{(X,W)})$.
Identifying $\mathbb{P}^N$ with $\mathbb{P}(W)$ (the set of
codimension $1$ vector subspaces of $W$) we see that
$\text{Osc}^k_x(X)$ is a linear subspace of $\mathbb{P}^N$. To
avoid that it fills up the whole ambient space we assume that $N$
is large enough. For instance, to discuss osculation for surfaces,
i.\ e., $k=n=2$, a reasonable assumption is that $N \geq 6$ or
even $5$, depending on the regularity of the surface we are
dealing with. Recalling that $\text{rk}(J_kL)=\binom{k+n}n -1$ we
have $\dim \big( \text{Osc}^k_x(X) \big) \leq \text{min}\{ N,
\binom{k+n}n -1\}$. Let $\mathcal{U} \subseteq X$ be the Zariski
dense open subset where the rank of the homomorphism
$j_{k,x}^{(X,W)}: W \to (J_kL)_x$ attains its maximum, say
$s(k)+1$. The $k$-th inflectional locus of $(X,W)$ is defined by
$\Phi_k(X)=X \setminus \mathcal{U}$. So $x \in \Phi_k(X)$ if and
only if $\dim \big( \text{Osc}^k_x(X) \big) < s(k)$. By
{\it{flex}} we simply mean a point in $\Phi_2(X)$, while a {\it
higher flex} is a point of $\Phi_k(X)$ with $k > 2$. We say that
$X$ is {\it uninflected} to mean that $\Phi_2(X)=\emptyset$. Of
course $\Phi_h(X) \subseteq \Phi_k(X)$ for $h \leq k$. Let $n=1$.
If $N < k$, then clearly $\Phi_k(X)=X$. However, if $N \geq k$
then $\Phi_k(X) \subsetneq X$ (e.\ g., see \cite[p.\ 37, Ex
C-2]{bib1}). In particular, $\Phi_N(X)=\emptyset$ if and only if
$X$ is a rational normal curve \cite[p.\ 39, Ex C-14]{bib1}.
\par
Now let $x \in \mathcal{U}$. A hyperplane $H \in \mathbb{P}^{N
\vee}$ is said to be $k$-{\it{th osculating to}} $X$ {\it{at}} $x$
if $H \supseteq \text{Osc}_x^k(X)$. Then the $k$-{\it{th dual
variety}} $X_k^{\vee}$ of $(X,W)$ is defined as the closure in
$\mathbb{P}^{N \vee}$ of the locus parameterizing all $k$-th
osculating hyperplanes to $X$ at points of $\mathcal{U}$.
\par
By {\it scroll} we mean an embedded smooth projective variety $Y
\subset \mathbb{P}^N$ of dimension $n \geq 1$ endowed with a
morphism $\pi:Y \to C$ over a smooth curve $C$ such that $(f,
\mathcal{O}_{\mathbb{P}^N}(1)|_f) = (\mathbb{P}^{n-1},
\mathcal{O}_{\mathbb{P}^{n-1}}(1))$ for every fibre $f$ of $\pi$,
or the corresponding pair $(X,W)$ with $|W|$ very ample, such that
$Y=\varphi_W(X)$. Of course $X=C$ if $n=1$. We need to fix some
more notation.
\par
Let $(X,W)$ be a scroll. As is known, for any $k\geq 2$ we have a
strict inequality $\dim \big( \text{Osc}^k_x(X) \big)<
\binom{k+n}n -1$ at every point $x \in X$. In fact, there are
local coordinates $(u,v_2, \dots , v_n)$ around every point $x \in
X$ such that the homogeneous coordinates $x_i$ $(i=0, \dots ,N)$
of the points of the variety locally can be written as $x_i =
a_i(u)+ \sum_{j=2}^n v_j b_{ij}(u)$, where $a_i$ and $b_{ij}$ are
holomorphic functions of $u$. Since every section $\sigma \in W$
is a linear combination $\sigma = \sum_{i=0}^N \lambda_ix_i$ we
thus see that the second derivatives $\sigma_{v_jv_h}$ vanish at
every point. Then $\dim \big( \text{Osc}^2_x(X) \big)\leq 2n$, and
differentiating further up to the order $k$ we see that
\[\dim \big( \text{Osc}^k_x(X) \big)\leq nk \quad \text{for every $x \in X$}
.\]
\par

Finally, we set
$\mathbb{F}_e=\mathbb{P}(\mathcal{O}_{\mathbb{P}^1} \oplus
\mathcal{O}_{\mathbb{P}^1}(-e))$ to denote the Segre--Hirzebruch
surface of invariant $e$ ($e \geq 0$). Then, as in \cite[p.\
372]{bib5}, $C_0$ stands for a section of minimal
self-intersection and $f$ for a fibre.

\section{Decomposable scrolls and their flexes}

The situation we consider for the most part of this paper is
inspired by that in \cite{bib9} and \cite[Sec.\ 2]{bib7}. Let $C$
be a smooth curve of genus $g$. For $i=1, \dots, n$ let $\mathcal
{L}_i$ be a very ample line bundle on $C$ and let $V_i \subseteq
H^0(C,\mathcal{L}_i)$ be a vector subspace such that $|V_i|$ gives
rise to an embedding
\[\varphi_i: C \to \mathbb{P}^{r_i}=\mathbb{P}(V_i).\]
Set $C_i=\varphi_i(C)$. Let $V=\oplus_{i=1}^n V_i$, $\mathcal{E} =
\oplus_{i=1}^n \mathcal{L}_i$ and consider the projective bundle
$P=\mathbb{P}(\mathcal{E})$. By identifying $V$ with a vector
subspace of $H^0(P,L)$, where $L$ is the tautological line bundle
on $P$, we get an embedding
\[\varphi: P \to \mathbb{P}^N = \mathbb{P}(V).\]
We set $X=\varphi(P)$. According to \cite[p.\ 151]{bib7} we say
that $X$ is the {\it{decomposable scroll}} generated by $C_1,
\dots, C_n$. For a point $p \in C$, let $p_i=\varphi_i(p)\in C_i$.
Geometrically, $X$ is generated by the linear spaces $f_p:=\langle
p_1, \dots , p_n \rangle \cong \mathbb{P}^{n-1}$ as the point $p$
varies on $C$; note that all the linear spans $\langle C_i \rangle
= \mathbb{P}^{r_i}$ of the $C_i$'s are skew each other and
generate the whole ambient space $\mathbb{P}^N$. Let $t$ be a
local parameter on $C$ such that $\varphi_i(p) = (x_0(0), \dots ,
x_{r_i}(0))$ corresponds to $t=0$. Locally, around $p$, the
homomorphism $j_k^{C_i}:V_i \to J_k\mathcal{L}_i$ is represented
by the matrix
\[ M^i_k(t) = \begin{pmatrix} x_0(t) & x_1(t) & \dots & x_{r_i}(t)\\
x^{\prime}_0(t) & x^{\prime}_1(t) & \dots & x^{\prime}_{r_i}(t)\\
. & . & \dots & .\\
. & . & \dots & .\\
x^{(k)}_0(t) & x^{(k)}_1(t) & \dots & x^{(k)}_{r_i}(t)\\
\end{pmatrix} . \]
The linear space spanned by the row vectors of the matrix
$M^i_k(0)$ defines the $k$-th osculating space to $C_i$ at $p_i$.
Note that, if $k > r_i$, then every $k$-th osculating space to
$C_i$ is the whole space $\mathbb{P}^{r_i}=\mathbb {P}(V_i)$. Now
let $\lambda_1, \dots ,\lambda_n$ denote homogeneous coordinates
corresponding to a local trivialization of $\mathcal{E}$ around
$p$, and, for $\lambda_n\not=0$, set $v_i=\lambda_i/\lambda_n$.
Then $(t,v_1, \dots , v_{n-1})$ provide local coordinates on $X$
at a point $x \in f_p \setminus \langle p_1, \dots , p_{n-1}
\rangle$. Writing down the parametric equations for $X$ around
$f_p$ we can easily get the matrix $M^X_k(t,v_1, \dots ,v_{n-1})$
representing $j^X_k:V \to J_kL$ near $x$. Set
\[ \overline{M^i_{k-1}} = M^i_{\text{min}\{k-1, r_i\}}.\]

\begin{lemma}\cite[p.\ 152]{bib7} We have
\[ M^X_k(t,v_1, \dots ,v_{n-1})= \begin{pmatrix}
v_1M^1_k & v_2 M^2_k
& \dots & v_{n-1} M^{n-1}_k & M^n_k \\
{\overline{M^1_{k-1}}} & 0  & \dots & 0 & 0\\
0 & {\overline{M^2_{k-1}}}  & \dots & 0 & 0\\
. & . & \dots & . & .\\
. & . & \dots & . & .\\
0 & 0 & \dots & {\overline{M^{n-1}_{k-1}}} & 0\\
0 & 0 & \dots & 0 & 0
\end{pmatrix} (t) \quad. \]
\end{lemma}
\par
Let $x \in f_p \setminus \langle p_1, \dots , \widehat{p_s}, \dots
, p_n \rangle$, where $\widehat{\ } $  denotes suppression. Up to
reordering the $C_i$'s, there is no restriction if we suppose that
$s=n$, hence the matrix representing $j_k^X$ is that given by
Lemma 1.1. Sometimes, however, it is convenient to order the
$C_i$'s according to some criterion (e.\ g., in such a way that
$r_1 \leq r_2 \leq \dots \leq r_n$). In this case, we can write
$x=u_1p_1 + \dots + u_{s-1}p_{s-1} + p_s + u_s p_{s+1} + \dots +
u_{n-1} p_n$, and then, with respect to the local coordinates
$(t,u_1, \dots , u_{n-1})$ the matrix representing $j_k^X$ near
$x$ is the following
\[ \begin{pmatrix} u_1M^1_k & \dots & u_{s-1} M^{s-1}_k &
M^s_k & u_s M^{s+1}_k & \dots
& u_{n-1} M^n_k\\
{\overline{M^1_{k-1}}} & \dots & 0 & 0 & 0 & \dots & 0\\
. & \dots & . & . & . & \dots & .\\
. & \dots & . & . & . & \dots & .\\
0 & \dots & {\overline{M^{s-1}_{k-1}}} & 0 & 0 & \dots & 0\\
0 & \dots & 0 & 0 & 0 & \dots & 0\\
0 & \dots & 0 & 0 & {\overline{M^{s+1}_{k-1}}} & \dots & 0\\
. & \dots & . & . & . & \dots & .\\
. & \dots & . & . & . & \dots & .\\
0 & \dots & 0 & 0 & 0 & \dots & {\overline{M^n_{k-1}}}\\
\end{pmatrix} (t) \quad. \]
\par
We say that two matrices $A$ and $B$ of type $m \times n$ are row
equivalent if the vector subspace of $\mathbb{C}^n$ spanned by the
rows of $A$ is the same as that spanned by the rows of $B$.
\par
Here is an immediate application.

\begin{theorem} Let $X$ be a decomposable scroll generated by $C_1,
\dots, C_n$ and let $\Phi_2(X)$ be its inflectional locus.
\par \noindent $(1)$ The following three conditions are equivalent:
\begin{enumerate}
\item[{\rm(i)}] $\big(f_p \setminus \bigcup_{i=1}^n \langle p_1,
\dots, \widehat{p_i}, \dots , p_n \rangle \big) \cap
\Phi_2(X)\not= \emptyset$; \item[{\rm(ii)}] $p_i$ is a flex of
$C_i$ for every $i=1, \dots, n$; \item[{\rm(iii)}] $f_p \subseteq
\Phi_2(X)$.
\end{enumerate}
\par \noindent $(2)$ $p_i \in \Phi_2(X)$ if and only if it is a flex of
$C_i$. \par \noindent $(3)$ Let $x \in \Phi_2(X)$: if $x \in f_p
\setminus \langle p_1, \dots, \widehat{p_s}, \dots , p_n \rangle$
then $p_s$ is a flex of $C_s$.
\end{theorem}

\begin{proof} To prove (1) it is enough to show that (i)
$\Rightarrow$ (ii) $\Rightarrow$ (iii). Let $x \in f_p \setminus
\langle p_1, \dots , p_{n-1} \rangle$, so that we can write $x =
v_1p_1+ \dots v_{n-1}p_{n-1}+p_n$. Then $x \in \Phi_2(X)$ if and
only if $j^X_{2,x}:V \to (J_2L)_x$ has rank $< 2n+1$. Note that
$\text{rk}\big( \overline{M_1^i}(t) \big)=2$ and
\begin{equation} \label{1.2.1}
\text{rk}\big( M_2^i(t) \big) \geq 2  \end{equation} for every $i$
and for every $t$. Then Lemma 1.1 shows that $\text{rk}\big(
M_2^X(0,v_1, \dots , v_{n-1}) \big) < 2n+1$ if and only if both
$\text{rk}\big( M_2^n(0) \big)=2$ and
\[\text{rk}\big( (v_1M_2^1 \quad v_2M_2^2 \quad \dots \quad
v_{n-1}M_2^{n-1})(0) \big) = 2.\] The former condition says that
$j_{2,p}:V_n \to (J_2\mathcal{L}_n)_p$ has rank $2$, while by
(\ref{1.2.1}) the latter one is equivalent to saying that either
$v_i=0$ or $\text{rk}\big( M_2^i(0) \big)=2$ for every $i=1, \dots
,n-1$. In conclusion we have that $j_{2,p}:V_i \to (J_2
\mathcal{L}_i)_p$ has rank $2$ for $i=n$ and for every $i$ such
that $v_i\not=0$. So, if $x \in \Phi_2(X)$ is a general point as
in (i), we get (ii). On the other hand, if (ii) holds, then we see
that $f_p \setminus \langle p_1, \dots , p_{n-1} \rangle$, hence
its closure $f_p$, lies in $\Phi_2(X)$. So (1) is proved.
Moreover, the above argument proves the ``only if'' part of (2)
when $x=p_n$, and (3) in the special case $s=n$. As to the ``if''
part of (2), note that if $x=p_n$ then the matrix $M_2^X(0,0,
\dots,0)$ of Lemma 1.1 has the following special form:
\[ \begin{pmatrix}
0 & 0
& \dots & 0 & M^n_2 \\
{\overline{M^1_1}} & 0  & \dots & 0 & 0\\
0 & {\overline{M^2_1}}  & \dots & 0 & 0\\
. & . & \dots & . & .\\
. & . & \dots & . & .\\
0 & 0 & \dots & {\overline{M^{n-1}_1}} & 0\\
0 & 0 & \dots & 0 & 0
\end{pmatrix} \quad. \]
So, if $p_n$ is a flex of $C_n$ we get $\text{rk}\big( M_2^X(0,0,
\dots,0) \big)=2n$, since $M_2^n(0)$ has rank $2$. Now, let $x$ be
as in (3); so we can write $x=u_1p_1+ \dots u_{s-1}p_{s-1}+p_s+
u_sp_{s+1}+ \dots +u_{n-1}p_n$. Then one can easily see that the
matrix representing $j_2^X:V \to J_2L$ near $x$ is
\[ \begin{pmatrix}
u_1M^1_2 & \dots & u_{s-1} M^{s-1}_2 & M^s_2 & u_s M^{s+1}_2 &
\dots
& u_{n-1} M^n_2\\
{\overline{M^1_1}} & \dots & 0 & 0 & 0 & \dots & 0\\
. & \dots & . & . & . & \dots & .\\
. & \dots & . & . & . & \dots & .\\
0 & \dots & {\overline{M^{s-1}_1}} & 0 & 0 & \dots & 0\\
0 & \dots & 0 & 0 & 0 & \dots & 0\\
0 & \dots & 0 & 0 & {\overline{M^{s+1}_1}} & \dots & 0\\
. & \dots & . & . & . & \dots & .\\
. & \dots & . & . & . & \dots & .\\
0 & \dots & 0 & 0 & 0 & \dots & {\overline{M^n_1}}\\
\end{pmatrix} (t) \quad. \]
Thus the same argument as above works and shows that since $x \in
\Phi_2(X)$, $p_s$ must be a flex of $C_s$. This completes the
proof of (3) and (2).
\end{proof}

\begin{corollary} $X$ is uninflected if and only if $C_1, \dots, C_n$ are
uninflected.
\end{corollary}
\par
The same argument proving Theorem 1.2 says more.

\begin{proposition} For any $x \in \Phi_2(X)$ we have
\[\text{\rm{Osc}}^2_{x}(X) = \langle \text{\rm{Osc}}^1_{p_1}(C_1),
\dots , \text{\rm{Osc}}^2_{p_s}(C_s), \dots ,
\text{\rm{Osc}}^{1}_{p_n}(C_n)\rangle\] for some $s$, where, $p_s
\in \Phi_2(C_s)$. Moreover, $\text{\rm{Osc}}^2_{x}(X)$ is the same
linear $\mathbb{P}^{2n-1}$ for all $x \in \Phi_2(X) \cap f_p$.
\end{proposition}

\begin{proof} First, suppose that $x \not\in \langle p_1, \dots,
p_{n-1} \rangle$. Then $x=v_1p_1+ \dots +v_{n-1}p_{n-1}+p_n$. As
$x \in \Phi_2(X)$, the first block of rows in the matrix $M^X_2(0,
v_1, \dots ,v_{n-1})$ appearing in Lemma 1.1 for $k=2$ has rank
$2$. In particular, $p_n \in \Phi_2(C_n)$ by Theorem 1.2(3).
Moreover, either $v_i=0$ or $v_iM^i_2$ is row equivalent to
$\overline{M^i_1}$. Hence $M^X_2$ is row equivalent to the matrix
\[ \begin{pmatrix} 0 & 0
& \dots & 0 & M^n_2 \\
{\overline{M^1_1}} & 0  & \dots & 0 & 0\\
0 & {\overline{M^2_1}}  & \dots & 0 & 0\\
. & . & \dots & . & .\\
. & . & \dots & . & .\\
0 & 0 & \dots & {\overline{M^{n-1}_1}} & 0\\
0 & 0 & \dots & 0 & 0
\end{pmatrix}  \quad. \]
This means exactly that
\[\text{Osc}^2_{x}(X) = \langle \text{Osc}^1_{p_1}(C_1),
\dots , \text{Osc}^{1}_{p_{n-1}}(C_{n-1}),
\text{Osc}^{2}_{p_n}(C_n)\rangle.\] Next, suppose that $x \in
\langle p_1, \dots, p_{n-1} \rangle \setminus \langle p_1, \dots,
p_{n-2} \rangle$. Then, $x=v_1p_1+ \dots +v_{n-2}p_{n-2}+p_{n-1}$
can also be written as $x=u_1p_1 + \dots +u_{s-1}p_{s-1} + p_s+
u_sp_{s+1} + \dots +u_{n-1}p_n$, as done after Lemma 1.1, with
$s=n-1$ and $u_{n-1}=0$. Then look at the matrix appearing after
Lemma 1.1 in the present situation:
\[ \begin{pmatrix}
u_1M^1_2 & \dots & u_{n-2} M^{n-2}_2 & M^{n-1}_2 & 0 &\\
{\overline{M^1_1}} & \dots & 0 & 0 & 0\\
. & \dots & . & . & .\\
. & \dots & . & . & .\\
0 & \dots & {\overline{M^{n-2}_1}} & 0 & 0\\
0 & \dots & 0 & 0 & {\overline{M^n_1}}\\
\end{pmatrix} (0) \quad. \]
If $x \in \Phi_2(X)$, arguing as before we see that $p_{n-1}\in
\Phi_2(C_{n-1})$ and this matrix is row equivalent to
\[ \begin{pmatrix}
0 & \dots & 0 & M^{n-1}_2 & 0 &\\
{\overline{M^1_1}} & \dots & 0 & 0 & 0\\
. & \dots & . & . & .\\
. & \dots & . & . & .\\
0 & \dots & {\overline{M^{n-2}_1}} & 0 & 0\\
0 & \dots & 0 & 0 & {\overline{M^n_1}}\\
\end{pmatrix} (0) \quad. \]
This means that
\[\text{Osc}^2_{x}(X) = \langle \text{Osc}^1_{p_1}(C_1),
\dots , \text{Osc}^{1}_{p_{n-2}}(C_{n-2}),
\text{Osc}^{2}_{p_{n-1}}(C_{n-1}),
\text{Osc}^{1}_{p_n}(C_n)\rangle.\] Now, let $s \leq n-2$. By
repeating the argument for $x \in \langle p_1, \dots, p_s \rangle
\setminus \langle p_1, \dots, p_{s-1} \rangle$, we see that
$p_s\in \Phi_2(C_s)$ and $\text{Osc}^2_{x}(X)$ is the linear span
of $\text{Osc}^{2}_{p_s}(C_s)$ and the spaces
$\text{Osc}^1_{p_i}(C_i)$ for $i \not=s$. This proves the first
assertion. Now, note that all $\text{Osc}^{1}_{p_i}(C_i)$ are
lines. Moreover, as we have shown, $p_s$ is a flex for $C_s$,
hence $\text{Osc}^{2}_{p_s}(C_s)$ is also a line. Thus, for any $x
\in X$, $\text{Osc}^2_{x}(X)$ is the linear space generated by the
$n$ tangent lines to $C_i$ at $p_i$ for $i=1, \dots ,n$. Note that
they generate a $\mathbb{P}^{2n-1}$. It turns out that
$\text{Osc}^2_{x}(X)$ is the same $\mathbb{P}^{2n-1}$ for all $x
\in \Phi_2(X) \cap f_p$.
\end{proof}

\section{Higher flexes and fibres}

Let $X$ be a decomposable scroll over a smooth curve $C$ generated
by $C_1, \dots , C_n$ as in Section 1, and let $f_p=\langle p_1,
\dots , p_n \rangle$ be the fibre over $p \in C$. In this section
we explore some connections between the higher inflectional loci
$\Phi_k(X)$ and the fibres of $X$.

\begin{remark} We have
\begin{equation} \label{2.1.1}
\text{Osc}^k_{p_s}(X) = \langle \text{Osc}^{k-1}_{p_1}(C_1), \dots
, \text{Osc}^k_{p_s}(C_s), \dots ,
\text{Osc}^{k-1}_{p_n}(C_n)\rangle \end{equation} for any $s=1,
\dots , n$ (the only $k$-th osculating space on the right hand is
that at $p_s$). In particular, if $p_s \in \Phi_k(C_s)$, then $p_s
\in \Phi_k(X)$.
\end{remark}

\begin{proof} Up to reordering the curves we can suppose that $s=n$.
Then the matrix representing $j_{k, p_s}^X$ is, according to Lemma
1.1,
\begin{equation} \label{2.1.2}
\begin{pmatrix}
0 & \dots & 0 & M^n_k \\
{\overline{M^1_{k-1}}} & \dots & 0  & 0\\
. & \dots & . & .\\
. & \dots & . & .\\
0 & \dots & {\overline{M^{n-1}_{k-1}}} & 0\\
0 & \dots & 0 & 0
\end{pmatrix} \quad.
\end{equation}
This proves the first assertion. Note that all linear spaces
appearing on the right hand of (\ref{2.1.1}) are skew each other.

Then the second assertion follows from the inequality:

\begin{equation}
\begin{split}
\dim\big( \text{Osc}^k_{p_s}(X) \big)
    &\leq (n-1)(k-1) + \dim\big( \text{Osc}^k_{p_s}(C_s) \big) + (n-1) \\
    & < (n-1)k + k = nk.  \end{split}
\end{equation}
\end{proof}

As to the converse, if $p_s \in \Phi_k(X)$, we cannot claim that
$p_s \in \Phi_k(C_s)$ if $k > 2$. However, we have

\begin{remark} If $p_s \in \Phi_k(X)$, then either $p_s \in
\Phi_k(C_s)$, or $p_j \in \Phi_{k-1}(C_j)$ for some $j \not= s$.
\end{remark}

\begin{proposition} Let $p_i \in \Phi_k(C_i)$ for $i=1, \dots ,
\widehat{s}, \dots ,n$, where $\widehat{\ }$ denotes suppression.
Then for every point $x \in f_p \setminus \langle p_1, \dots ,
\widehat{p_s}, \dots , p_n \rangle$ we have
\[\text{\rm{Osc}}^k_x(X) = \langle \text{\rm{Osc}}^{k-1}_{p_1}(C_1),
\dots , \text{\rm{Osc}}^k_{p_s}(C_s), \dots ,
\text{\rm{Osc}}^{k-1}_{p_n}(C_n)\rangle \] (the only $k$-th
osculating space on the right hand is that at $p_s$).
\end{proposition}

\begin{proof} Up to reordering we can suppose that $s=n$. Due to the
assumption, we have
$\text{Osc}^k_{p_i}(C_i)=\text{Osc}^{k-1}_{p_i}(C_i)$ for $i=1,
\dots, n-1$. This means that the two matrices $M^i_k$ and
$\overline{M^i_{k-1}}$ are row equivalent for $i=1, \dots , n-1$.
Now look at the matrix $M$ of Lemma 1.1. By subtracting suitable
linear combinations of the subsequent rows from the first block of
rows we see that $M$ is row equivalent to the matrix
(\ref{2.1.2}). This proves the assertion.
\end{proof}
\par
The same argument proves more. Actually, assume that $p_{i_j} \in
\Phi_k(C_{i_j})$ for $j=1, \dots ,s$ and set $\Lambda = \langle
p_{i_1}, \dots , p_{i_s} \rangle$. Up to reordering we can suppose
that $(i_1, \dots ,i_{s-1},i_s) = (1, \dots , s-1, n)$. Then for
any $x \in \Lambda \setminus \langle p_1, \dots ,p_{s-1} \rangle$
we can write $x = v_1p_1 + \dots v_{s-1}p_{s-1}+p_n$. Arguing as
in the proof of Proposition 2.3 we have that the matrices $M^i_k$
and $\overline{M^i_{k-1}}$ are row equivalent for $i=1, \dots ,
s-1$. Now look at the matrix $M$ of Lemma 1.1, representing
$j_{k,x}$. The first block of rows of $M$ is
\[ \begin{pmatrix} v_1M^1_k & \dots & v_{s-1} M^{s-1}_k & 0 & \dots &
0 & M^n_k
\end{pmatrix} . \]
By subtracting suitable linear combinations of the subsequent rows
of $M$ from the first block we see that $M$ is row equivalent to
the matrix in (\ref{2.1.2}). Now, since also $p_n \in \Phi_k(C_n)$
we have $\text{rk}\big( M^n_k \big) < k+1$ and then the same
computation done to prove Remark 2.1 holds at $x$, giving
$\dim\big( \text{Osc}^k_x(X) \big)<nk$. Thus $\Lambda \setminus
\langle p_1, \dots ,p_{s-1} \rangle \subseteq \Phi_k(X)$. On the
other hand $\Phi_k(X) \cap f_p$ is a Zariski closed subset, hence
$\Lambda \subseteq \Phi_k(X)$.

\par
Now suppose that $(i_1, \dots , i_s) = (1, \dots , s)$, with $s
\leq n-1$ and $p_n \not\in \Phi_k(C_n)$. Then $\text{rk}\big(
M^n_k \big)= k+1$, and the same argument as above applied to any
point $x \in \langle p_1, \dots p_s, p_n \rangle \setminus
\Lambda$ shows that
\[\dim \big(\text{Osc}^k_x(X) \big) = \sum_{i=1}^{n-1} \text{rk}
\big( \overline{M^i_{k-1}}\big) + (k+1) -1.\] In particular, if
$p_i \not\in \Phi_{k-1}(C_i)$ for $i=1, \dots, s$, then all the
first $n-1$ summands are equal to $k$, hence $\dim \big
(\text{Osc}^k_x(X) \big) = nk$, and so $x \not\in \Phi_k(X)$. This
proves the following

\begin{proposition} If $p_{i_j} \in \Phi_k(C_{i_j})$ for $j =1, \dots
,s \leq n$, then $\langle p_{i_1}, \dots , p_{i_s} \rangle
\subseteq \Phi_k(X)$. Moreover, if $p_{i_j} \in \Phi_k(C_{i_j})
\setminus \Phi_{k-1}(C_{i_j})$ for $j=1, \dots , s \leq n-1$ and
$p_{i_j} \not\in \Phi_k(C_{i_j})$ for $j=s+1, \dots, n$, then
$\Phi_k(X) \cap f_p = \langle p_{i_1}, \dots ,p_{i_s} \rangle$.
\end{proposition}

\begin{corollary} {\rm{i)}} If $p_i \in \Phi_k(C_i)$ for every $i =1, \dots ,n$,
then $f_p \subseteq \Phi_k(X)$. \par {\rm{ii)}} If $f_p \subseteq
\Phi_k(X)$, then $p_i \in \Phi_k(C_i)$ for some $i$.
\end{corollary}

\begin{proof} i) is obvious; ii) follows from Remark 2.2, taking into
account the inclusion $\Phi_{k-1}(C_j) \subseteq \Phi_k(C_j)$
\end{proof}

\par In particular, if $f_p \subseteq \Phi_k(X)$ and $k > 2$, we see that not
necessarily $p_i \in \Phi_k(C_i)$ for all $i$'s. For $n=2$ we can
be more explicit.

\begin{proposition} Let $n=2$ and $k \geq 2$. Then $f_p \subseteq \Phi_k(X)$
if and only if either
\begin{enumerate}
\item[{\rm a)}] $p_i \in \Phi_{k-1}(C_i)$ for some $i$, or
\item[{\rm b)}] $p_i \in \Phi_k(C_i)$ for $i=1,2$.
\end{enumerate}
\end{proposition}

\begin{proof} Let $f_p \subseteq \Phi_k(X)$. By Corollary 2.5, ii),
up to reordering, we can suppose that $p_1 \in \Phi_k(C_1)$. Then,
by Proposition 2.3, for every $x \in f_p \setminus \{p_1\}$, we
have

\begin{equation} \label{2.6.1}
\text{\rm{Osc}}^k_x(X) = \langle \text{\rm{Osc}}^{k-1}_{p_1}(C_1),
\text{\rm{Osc}}^k_{p_2}(C_2)\rangle . \end{equation} Hence

\begin{equation} \label{2.6.2}
\dim\big( \text{\rm{Osc}}^k_x(X) \big) = \dim \big(
\text{\rm{Osc}}^{k-1}_{p_1}(C_1) \big)+ \dim \big(
\text{\rm{Osc}}^k_{p_2}(C_2) \big)+1 .  \end{equation} Since $x
\in \Phi_k(X)$ this shows that either $p_1 \in \Phi_{k-1}(C_1)$,
case a), or $p_2 \in \Phi_k(C_2)$, case b). To prove the converse,
in both cases a) and b), up to renaming, we can assume that $p_1
\in \Phi_k(C_1)$. Then Proposition 2.3 gives again (\ref{2.6.1})
for any $x \in f_p \setminus \{p_1\}$ and then (\ref{2.6.2}) shows
that
\[
\dim \big( \text{\rm{Osc}}^k_x(X) \big) \leq \begin{cases}
k-2+k+1 \quad \text{in case a)},\\
k-1+(k-1)+1 \quad \text{in case b).}\\
\end{cases} \]
Hence $f_p \setminus \{p_1\} \subseteq \Phi_k(X)$ in both cases,
and then, taking the closure, we get $f_p \subseteq \Phi_k(X)$.
\end{proof}

\begin{theorem} Let $n=2$. Suppose that $x \in \Phi_k(X)$ and let
$f_p$ be the fibre of $X$ containing $x$.

\begin{enumerate}
\item[{\rm i)}] If $x \not= p_1, p_2$, then $f_p \subseteq
\Phi_k(X)$; \item[{\rm ii)}] if $x=p_i$, then either $f_p
\subseteq \Phi_k(X)$ or $p_i \in \Phi_k(C_i)$.
\end{enumerate}
\end{theorem}

\begin{proof} If $x \not= p_1$, then we can write $x = vp_1+p_2$.
According to Lemma 1.1, $j_{k,x}^X$ is represented by the
following matrix
\[ M = \begin{pmatrix}
vM^1_k & M^2_k \\
{\overline{M^1_{k-1}}} & 0
\end{pmatrix} \quad. \]
Since $x \in \Phi_k(X)$, $M$ has rank $\leq 2k$. This implies
either
\begin{enumerate}
\item[$\alpha$)] $\text{rk}\big( {\overline{M^1_{k-1}}}\big) <
k-1$, i.\ e., $p_1 \in \Phi_{k-1}(C_1)$, or \item[$\beta$)]
$\text{rk} \big( vM^1_k \quad M^2_k \big) < k$.
\end{enumerate}
Condition $\beta$) in turn implies both $\text{rk}\big( M^1_k
\big) < k$ and $\text{rk} \big( M^2_k  \big) < k$. The latter
condition means that $p_2 \in \Phi_k(C_2)$, while the former is
equivalent to saying that
\[ \text{either} \quad v = 0, \qquad \text{or} \qquad \text{rk}
\big(M^1_k \big) < k.\] In other words, either $x=p_2$ or $p_1 \in
\Phi_k(C_1)$. In conclusion, if $x \not= p_1, p_2$ then either
\begin{enumerate}
\item[$\alpha$)] $p_1 \in \Phi_{k-1}(C_1)$, or \item[$\beta$)]
$p_i \in \Phi_k(C_i)$ for $i=1,2$.
\end{enumerate}
In both cases $f_p \subseteq \Phi_k(X)$ by Proposition 2.6. This
proves i). Now let $x = p_i$. By Remark 2.2 either $p_i \in
\Phi_k(C_i)$ or $p_j \in \Phi_{k-1}(C_j)$ for $j \not=i$. But in
the latter case Proposition 2.6 says that $f_p \subseteq
\Phi_k(X)$ again. This proves ii).
\end{proof}

\begin{corollary} Let $n=2$. Then $\Phi_k(X)= \emptyset$ if and only
if $\Phi_k(C_i) = \emptyset$ for $i=1,2$.
\end{corollary}

\begin{proof} If $p_i \in \Phi_k(C_i)$ for some $i$, we know that
$p_i \in \Phi_k(X)$ by Remark 2.2. This proves the ``only if
part''. To prove the ``if'' part suppose, by contradiction, that
$x \in \Phi_k(X)$, and let $f_p$ be the fibre of $X$ through $x$.
By Theorem 2.7 either $f_p \subseteq \Phi_k(X)$ or $p_i \in
\Phi_k(C_i)$ for some $i$. In both cases, taking into account
Proposition 2.6, we see that $\Phi_k(C_i) \not= \emptyset$ for
some $i$. But this is a contradiction.
\end{proof}

\par Now suppose that $r_1 \leq r_2$, where $\langle C_i
\rangle =\mathbb P^{r_i}$. If $r_1 < k-1$, then we have
$\dim(\text{Osc}^{k-1}_y(C_1)) \leq r_1 < k-1$ for every $y \in
C_1$. In other words, $\Phi_{k-1}(C_1)=C_1$ and therefore
$\Phi_k(X)=X$ by Proposition 2.6. Let $r_1 \geq k-1$. If $r_2 < k$
(i.\ e., $r_1=r_2=k-1$), then $\Phi_k(C_i)=C_i$ for $i=1,2$, hence
$\Phi_k(X)=X$ again, by Proposition 2.6. If $r_1=k-1$ but $r_2
\geq k$ then $\Phi_k(C_1)=C_1$ but $\Phi_k(C_2) \subsetneq C_2$
(e.\ g., see \cite[p.\ 37, Ex.\ C-2]{bib1}); so $\Phi_k(X)$
contains every fibre of $X$ passing through a point of either
$\Phi_{k-1}(C_1)$ or $\Phi_k(C_2)$. Taking into account also
Remark 2.1, we thus get the following

\begin{corollary} Let $n=2$ and suppose that $r_1 \leq r_2$. If
$\Phi_{k-1}(C_1)=\emptyset$ but $\Phi_k(C_1)=C_1$, then
$\Phi_k(X)$ consists of $C_1$ plus the fibres containing a point
of $\Phi_k(C_2)$.
\end{corollary}

\section{Surface scrolls: examples and applications}

In this section we focus on the surface case ($n=2$). Let $X
\subset \mathbb{P}^N$ be a decomposable surface scroll as in
Section 1. We present some examples concerned with the dimension
that $\text{Osc}_x^k(X)$ can have at some point $x$, and with the
structure of $\Phi_2(X)$, focusing in particular on the case of
non-normal rational scrolls. First it is useful to recall the
situation for normal rational scrolls.

\begin{example} Notation as in Section 1; let $C =\mathbb{P}^1$,
$\mathcal{E} = \mathcal{O}_{\mathbb{P}^1}(r_1) \oplus
\mathcal{O}_{\mathbb{P}^1}(r_2)$, with $1 \leq r_1 \leq r_2$, and
let $X \subset \mathbb{P}^N$ be the image of
$\mathbb{P}(\mathcal{E})$ in the embedding given by complete
linear system associated with the tautological line bundle $L$.
Note that $N=r_1+r_2+1$. Let $p = (t_0:t_1) \in \mathbb{P}^1$ and
set $t=t_1/t_0$ (or $t_0/t_1$). At any point $x \in X \setminus
C_1$ we can use local coordinates $(t,v)$ to write $x=vp_1+p_2$ on
the fibre $f_p$; then, according to Lemma 1.1, the homomorphism
$j_k^X:H^0(X,L) \to J_kL$ is represented near $x$ by the matrix
\[M_k^X(t,v) = \begin{pmatrix}
v M^1_k & M^2_k\\
{\overline{M^1_{k-1}}} & 0\\
0 & 0
\end{pmatrix} (t) \quad. \]
Note that
\[\text{rk}\big(M^2_k(t) \big)= \min\{k+1,r_2+1\}.\]
Moreover
\[\text{rk}\big(\overline{M^1_{k-1}}(t) \big)= \begin{cases}
\text{rk}(M^1_{k-1})=k \quad \text{if $k-1
\leq r_1$},\\
\text{rk}(M^1_{r_1})=r_1+1 \quad \text{otherwise.}\\
\end{cases} \]
It follows that
\begin{equation}
\begin{split}
\text{rk}(j^X_{k,x})&=\text{rk}\big(M_k^X(t,v) \big)\\
&=\text{rk}\big(M^2_k(t) \big)+
\text{rk}\big(\overline{M^1_{k-1}}(t) \big)\\
&=\min\{k+1,r_2+1\}+\min\{k, r_1+1\}. \end{split}
\end{equation}
Therefore
\begin{equation} \label{3.1.1}\dim \big( \text{Osc}^k_x(X) \big)=
\begin{cases} 2k \qquad \qquad \text{if $k
\leq r_1+1$},\\
k+r_1+1 \quad \text{if $r_1+1 \leq k \leq r_2$,}\\
r_1+r_2+1 \quad \text{if $k \geq r_2$.}\\
\end{cases} \end{equation}
Note that at any point $x \in X \setminus C_1$ the dimension of
the $k$-th osculating space can be strictly smaller than $2k$.
This is obvious when $N < 2k$, but it can happen also for $k \leq
\left[\frac{N-1}{2}\right]$, e.\ g., for a very unbalanced
rational normal scroll (i.\ e., with invariant $e:=r_2-r_1$ very
large). In fact, for $k \leq \left[\frac{N-1}{2}\right]$ we have
$\dim \big( \text{Osc}^k_x(X) \big) < 2k$ if $r_1+1 < k$ from
(\ref{3.1.1}). This means $N+1-e < 2k \leq N-1$, hence $e \geq 3$
is enough. For $k \geq 2$, (\ref{3.1.1}) also shows that
\[k+2 \leq \dim\big( \text{Osc}^k_x(X) \big) \leq 2k\]
at any point $x \in X \setminus C_1$. In particular, letting $k=2$
we see that $\dim \big( \text{Osc}^2_x(X) \big)=4$ for any $x \in
X \setminus C_1$.
\end{example}
\par We want to stress that $X$ was linearly
normal in the example above. Here is an enlightening example
showing how small the dimension of $\text{\rm{Osc}}^k_x(X)$ can be
at some point $x$, for any $k$, when we drop linear normality.

\begin{example} Fix integers $k \geq 2$ and $r =r_2 \geq 3$. Let
$C=\mathbb{P}^1$, $\mathcal{E} = \mathcal{L}_1 \oplus
\mathcal{L}_2$, where $\mathcal{L}_1 =
\mathcal{O}_{\mathbb{P}^1}(1)$, $\mathcal{L}_2 =
\mathcal{O}_{\mathbb{P}^1}(k+r-1)$, and let $V_1=H^0(\mathbb{P}^1,
\mathcal{L}_1) =\langle t_0, t_1 \rangle$,
\[V_2=\langle t_0^{k+r-1},\ t_0^{k+r-2}t_1,\
t_0^{r-2}t_1^{k+1},\ \dots \ ,t_0t_1^{k+r-2},\ t_1^{k+r-1}
\rangle.\] Note that $\varphi_2:\mathbb{P}^1 \to \mathbb{P}^r$
defines an embedding, which is not linearly normal, since $\dim
V_2 = r+1 < h^0(\mathcal{L}_2)$. Then $X \subset
\mathbb{P}^{r+2}$, defined as in Section 1, is a rational
non-normal scroll. Let $L$ be the hyperplane bundle and let $V
\subset H^0(X,L)$ be the subspace giving rise to the embedding.
Note that at the point $p \in \mathbb{P}^1$, corresponding to
$(t_0:t_1)=(1:0)$ we have
\[|V_2-2p|= \dots =|V_2-(k+1)p|.\] This
means that for every $h$, ($2 \leq h \leq k$) the homomorphism
\[j_{h,p}^{C_2}:V_2 \to (J_k \mathcal{L}_2)_p\] has a $2$-dimensional image
(isomorphic to $(J_1 \mathcal{L}_2)_p$), i.\ e.,
$\text{rk}(j_{k,p}^{C_2})=2$. On the other hand, at every point $q
\in \mathbb{P}^1$ it is obvious that $\text{rk}(j_{h,q}^{C_1})=2$
for any $h \geq 1$. Now, let $x \in f_p$. If $x \in f_p \setminus
\{p_2\}$, Proposition 2.3 shows that
\[\text{Osc}_x^h(X) = \langle \text{Osc}_{p_1}^h(C_1),
\text{Osc}_{p_2}^{h-1}(C_2)\rangle\] for any $h=2, \dots, k$. On
the other hand, Remark 2.1 tells us that \[\text{Osc}_{p_2}^h(X) =
\langle \text{Osc}_{p_1}^{h-1}(C_1),
\text{Osc}_{p_2}^h(C_2)\rangle\] for any $h \geq 2$. In both cases
$\text{Osc}_x^h(X)$ is the linear span of two skew lines, namely
$C_1$ and $\text{Osc}_{p_2}^1(C_2)$; hence
\[\text{Osc}^k_x(X)=\text{Osc}^{k-1}_x(X)= \dots = \text{Osc}^2_x(X) \quad \text{\rm{for every
$k\geq 3$}},\] at every point $x \in f_p$. In particular,
\[\dim \big( \text{Osc}^k_x(S) \big)=3 \quad \text{for all
$k \geq 2$ at any point $x \in f_p$}. \]
\end{example}

We recall that if $X \in \mathbb{P}^N$ is any scroll of dimension
$n$, then $\dim \big( \text{Osc}^2_x(X) \big) \geq n+1$
\cite{bib2} (see also \cite{bib6} for $n=2$).

\begin{example} Let $C=\mathbb{P}^1$, $\mathcal{L}_1 =
\mathcal{O}_{\mathbb{P}^1}(m)$ with $m \geq 2$, $V_1 =
H^0(\mathbb{P}^1, \mathcal{L}_1)$, and consider $\mathcal{L}_2 =
\mathcal{O}_{\mathbb{P}^1}(d)$ with $d \geq m+2$. The vector space
$H^0(\mathbb{P}^1, \mathcal{L}_2)$ defines an embedding of $C$ as
a rational normal curve $\Gamma \subset \mathbb{P}^d$. Projecting
$\Gamma$ from a general linear space $T$ of dimension $d-m-2$ to a
$\mathbb{P}^{m+1}$ we get an embedding. Let $V_2 = V(T)$ be the
vector subspace of $H^0(\mathbb{P}^1, \mathcal{L}_2)$
corresponding to this embedding. Let $C_i$ be the image of $C$ in
the embedding defined by $V_i$, $i=1,2$, and in the space
$\mathbb{P}^{2m+2}=\mathbb{P}(V_1 \oplus V_2)$ consider the
decomposable rational scroll $X$ generated by $C_1$ and $C_2$.
Note that $X = \mathbb{F}_{d-m}$. We claim that $\Phi_m(X)=
\emptyset$. Of course $\Phi_m(C_1) = \emptyset$. Let $O_{\Gamma}$
be the $m$-th osculating developable of $\Gamma$ (i.\ e., the
variety generated by the linear spaces $\text{Osc}^m_x(\Gamma)$,
as $x$ varies on $\Gamma$). Note that $\dim (O_{\Gamma}) = m+1$,
hence $T \cap O_{\Gamma}=\emptyset$ for a general $T$. Since no
osculating space $\text{Osc}^m_x(\Gamma)$ meets the center of
projection $T$, we conclude that $\Phi_m(C_2) = \emptyset$. Then
the claim follows from Corollary 2.8.
\end{example}

\par
Let us recall the following conjecture of Piene--Tai \cite{bib10}.
Let $S \subset \mathbb P^N$ ($N \geq 5$) be a non-degenerate
smooth projective surface such that $\dim \big(\text{Osc}^k_x(S)
\big) \leq 2k$ for all points $x \in S$ and for every $k$, with
equality for $k=[\frac{N-1}{2}]$, where $[\ ]$ stands for the
greatest integer function.
\begin{enumerate}
\item[(i)] If $N$ is odd, then $S$ is the balanced rational normal
scroll of degree $N-1$ (i.\ e., $S$ is $\mathbb{F}_0$ embedded by
$|C_0+[\frac{N-1}{2}]f|$). \item[(ii)] If $N$ is even, then $S$ is
the semibalanced rational normal scroll of degree $N-1$ (i.\ e.,
$S$ is $\mathbb{F}_1$ embedded by $|C_0+([\frac{N-1}{2}] + 1)f|$).
\end{enumerate}

Part (i) of this conjecture is true, as proved in \cite{bib3},
while part (ii) is not (see \cite[Theorem A and comment after
Corollary 2.3]{bib6}). Example 3.3 provides a new series of
counterexamples to the even dimensional part of the conjecture. We
want to stress that all these scrolls are decomposable, while
those appearing in \cite[Theorem A]{bib6} are not, all being
isomorphic to the elliptic $\mathbb{P}^1$-bundle of invariant
$-1$. Moreover, we have the following characterization, which
provides more information in order to correct the conjecture.

\begin{theorem} Let $X \subset \mathbb{P}^{2m+2}$ ($m \geq 2$) be a
decomposable scroll with $n=2$ such that $\Phi_m(X)=\emptyset$.
Then either $X$ is the semibalanced rational normal scroll of
degree $m+1$, or $X$ is of the type described in Example $3.3$.
\end{theorem}

\begin{proof} By Corollary 2.8 it must be $\Phi_m(C_i)= \emptyset$,
for $i=1,2$. In particular, $C_1$ cannot be a line, hence $r_1 =
\dim \big( \langle C_1 \rangle \big) \geq 2$. We can assume that
$r_1 \leq r_2$ and then from $r_1+r_2+1 = 2m+2$ we get that $r_1
\leq m$. As $\Phi_m(C_1)=\emptyset$, this implies that $r_1=m$,
and then $r_2= m+1$. So $C_1$ is a rational normal curve of degree
$m$ in $\mathbb{P}^m$ while $C_2$ is either the rational normal
curve of degree $m+1$ in $\mathbb{P}^{m+1}$ or any other rational
non-normal curve of some degree $d \geq m+2$ in
$\mathbb{P}^{m+1}$. In the former case $X$ is the semibalanced
rational normal scroll. In the latter, $C_2$ is obtained by
projecting a rational normal curve of degree $d$ in $\mathbb{P}^d$
to $\mathbb{P}^{m+1}$ from a general center as in Example 3.3.
\end{proof}
\par
The examples in the next part of this section are concerned with
$\Phi_2(X)$. First we would like to stress that for the cubic
scroll $X \subset \mathbb{P}^4$ the inflectional locus $\Phi_2(X)$
consists exactly of the generating line $C_1$. In fact this is the
only semi-balanced rational normal scroll which is not
uninflected. As to quartic rational normal scrolls in
$\mathbb{P}^5$ the situation is also well known \cite{bib11}. Let
us note that the one isomorphic to $\mathbb{F}_0$ is uninflected
according to Corollary 1.3, being generated by two conics $C_1$,
$C_2$. On the other hand, the one isomorphic to $\mathbb{F}_2$ is
generated by a line $C_1$ and a rational normal cubic $C_2$, which
has no flexes. Hence, according to Theorem 1.2(1) its inflectional
locus $\Phi_2$ consists exactly of $C_1$.
\par
\begin{example} We consider quintic non-normal rational scrolls in
$\mathbb{P}^5$. Let $X$ be as in Example 3.3, with $m=1$ and
$d=4$. According to Theorem 1.2(1), the inflectional locus
$\Phi_2(X)$ consists of the line $C_1$ and the fibres passing
through the flexes of the non-normal quartic rational curve $C_2$.
Now the center of projection $T$ is a point. If $T \not\in
O_{\Gamma}$, then $C_2$ has no flexes, as we said, and then
$\Phi_2(X) = C_1$. On the other hand, if $c \in O_{\Gamma}$, then
$C_2$ has $\epsilon$ flexes, where $\epsilon$ is the number of
osculating planes to $\Gamma$ passing through $T$. According to
the enumerative formula counting the weighted number of
$2$-osculating lines and $3$-osculating planes to $C_2$
\cite[Theorem 3.2]{bib8} we can see that $\epsilon = 1$ or $2$.
Depending on this, $\Phi_2(X)$ consists of $C_1$ plus one or two
fibres.
\end{example}
\par
\begin{example} In the same vein we can construct non-normal
rational scrolls having a finite inflectional locus. Let $C,
\mathcal{L}_2, V_2$ be as in the previous example but now put
$\mathcal{L}_1 =\mathcal{O}_{\mathbb{P}^1}(2)$ and
$V_1=H^0(\mathbb{P}^1, \mathcal{L}_1)$. Again let $C_i$ be the
image of $C$ in the embedding defined by $V_i$, $i=1,2$, and in
the space $\mathbb{P}^6=\mathbb{P}(V_1 \oplus V_2)$ consider the
decomposable sextic rational scroll $X$ generated by $C_1$ and
$C_2$. Now $C_1$ is a conic, hence it has no flexes. On the other
hand $C_2$ has $\epsilon=1$ or $2$ flexes provided that the
projection of $\Gamma$ giving rise to $C_2$ is made from a center
$T \in O_{\Gamma}$. Therefore, according to Theorem 1.2 ((1) and
(2)), $\Phi_2(X)$ consists of one or two points (the flexes of
$C_2$).
\end{example}
\par
\begin{example} Let $C$ be a smooth curve of genus $1$,
$\mathcal{L}_1=\mathcal O_{C}(3p)$, for some point $p \in C$, $V_1
= H^0(C, \mathcal{L}_1)$. Then $C_1=\varphi_1(C)$ is a smooth
plane cubic having exactly $9$ flexes, one of which is
$p_1:=\varphi_1(p)$. Now let $\mathcal{L}_2$ be a line bundle of
degree $5$ on $C$. The vector space $H^0(C, \mathcal{L}_2)$
defines an embedding of $C$ in $\mathbb{P}^4$ whose image, say
$\Gamma$, is a quintic normal elliptic curve: then $\Gamma$ has no
flexes (but $25$ hyperflexes). Projecting $\Gamma$ from a point $c
\in \mathbb{P}^4 \setminus \text{Sec}(\Gamma)$ to a $\mathbb{P}^3$
we get an embedding; let $V_2 = V(c)$ be the corresponding vector
subspace of $H^0(C, \mathcal{L}_2)$ and let $C_2$ be the image of
$C$ in the embedding $\varphi_2:C \to \mathbb{P}^3$ defined by
$V_2$. In the space $\mathbb{P}^6=\mathbb{P}(V_1 \oplus V_2)$
consider the decomposable elliptic scroll $X$ generated by $C_1$
and $C_2$. It can happen that $p_2:=\varphi_2(p)$ is a flex of
$C_2$ or not. According to Theorem 1.2(1), in the former case the
whole fibre $f_p$ is in the inflectional locus $\Phi_2(X)$, while
in the latter we have $\Phi_2(X) \cap f_p=\{p_1\}$. Note that if
$c$ is general enough, then $C_2$ has no flexes and therefore $X$
has only $9$ flexes: those of $C_1$.
\end{example}

\section{The second discriminant locus of decomposable scrolls}

Let $(X,L,W)$ be as at the beginning of Section 0, and let
$\mathcal{U} \subseteq X$ be the Zariski dense open subset of $X$
where $j_{k,x}^{(X,W)}: W \to (J_kL)_x$ attains the maximum rank
$s(k)+1$. If $x \in \mathcal{U}$, the fact that $H \in |W|$ is a
$k$-th osculating hyperplane to $X$ at $x$ is equivalent to the
fact that $H = (\sigma)_0$, where $\sigma \in W$ and
$j_{k,x}(\sigma)=0$. Equivalently, this means that $H \in
|W-(k+1)x|$, i.\ e., the hyperplane section cut out by $H$ on $X$
has a point of multiplicity $\geq (k+1)$ at $x$. Note however that
if $x \not\in \mathcal{U}$ and $H \in |W-(k+1)x)|$, this does not
necessarily mean that $H \in X_k^{\vee}$. Actually $H \in
X_k^{\vee}$ if and only if $H$ is a limit of $k$-th osculating
hyperplanes to $X$ at points of $\mathcal U$. On the other hand we
can consider the $k$-{\it{th discriminant locus}}
$\mathcal{D}_k(X,W)$ {\it{of}} $(X,W)$, which is defined as the
image of
\[\mathcal{J} :=\{(x,H)\in X \times |W|\ | \ H \in |W -(k+1)x|\}\]via
the second projection of $X \times |W|$. It parameterizes all
hyperplane sections of $X \subset \mathbb{P}^N$ admitting a
singular point of multiplicity $\geq k+1$; of course
$\mathcal{D}_k(X,W) \supseteq X_k^{\vee}$ with equality if and
only if $\mathcal{U} = X$, i.\ e., if and only if
\[ \dim \big( \text{Osc}_x^k(X) \big) = s(k)
\qquad \text{for every $x \in X$}.\] In general
$\mathcal{D}_k(X,W)$ contains some extra components coming from
the irreducible components of $\Phi_k(X)$.
\par
The discussion above says that $\mathcal{D}_k(X,W) = X_k^{\vee}$
if and only if $\Phi_k(X)=\emptyset$. From this point of view, the
characterization of balanced rational normal surface scrolls due
to Ballico, Piene and Tai \cite{bib3}, mentioned after Example
3.3, can be rephrased as follows.

\begin{proposition} Let $X \subset \mathbb{P}^N$ be any smooth surface,
where $N=2m+1 \geq 5$. Then $\mathcal{D}_m(X,W) = X_m^{\vee}$ if
and only if $X=\mathbb{P}^1 \times \mathbb{P}^1$ and
$W=H^0(\mathbb{P}^1 \times \mathbb{P}^1, \mathcal{O}_{\mathbb{P}^1
\times \mathbb{P}^1}(1,m))$.
\end{proposition}
\par
Now, let $X \subset \mathbb{P}^N=\mathbb{P}(V)$ be a decomposable
scroll as in Section 1. For simplicity we identify $X$ with the
corresponding abstract projective bundle $P$. So, we denote by
$\mathcal{D}_2(X,V)$ the second discriminant locus of $(P,L,V)$.
Its main component is the second dual variety $X_2^{\vee}$ of $X$.
Note that if $X$ is not linearly normal then $X_2^{\vee}$
corresponds to a suitable linear section of the second dual
variety of the linearly normal scroll giving rise to $X$ via the
projection to $\mathbb{P}^N$. Here, relying on the results of
Sections 1 and 2, we want to describe the extra components of
$\mathcal{D}_2(X,V)$. Of course we assume that $\Phi_2(X) \not=
\emptyset$. As a first thing we need to describe the irreducible
components of $\Phi_2(X)$.

\begin{proposition} Let $X \subset \mathbb{P}^N$ be a decomposable scroll
as in Section $1$, generated by $C_1, \dots ,C_n$, and assume that
$\Phi_2(X) \not= X$. Let $G$ be an irreducible component of
$\Phi_2(X)$. Then, up to reordering the curves $C_i$'s, either
\begin{enumerate}
\item[(1)] $G=\langle p_1, \dots ,p_s \rangle \subseteq f_p =
\langle p_1, \dots ,p_n \rangle$, or \item[(2)] $X$ is rational
and $G=C_1 \times \langle p_1, \dots ,p_s \rangle$ is the image of
$\mathbb{P}^1 \times \mathbb{P}^{s-1}$ via the Segre embedding.
\end{enumerate}
Moreover,
\[\text{\rm{Osc}}_x^2(X)=\langle \text{\rm{Osc}}_{p_1}^2(C_1),
\text{\rm{Osc}}_{p_2}^1(C_2), \dots, \text{\rm{Osc}}_{p_n}^1(C_n)
\rangle\] for all $x \in G$ in case $(1)$ and for all $x \in G
\cap f_p$ in case $(2)$. In particular, $\dim \big(
\text{\rm{Osc}}_x^2(X) \big)=2n-1$ for any $x \in G$ in both
cases.
\end{proposition}

\begin{proof} As $\Phi_2(X)\not=\emptyset$ it follows from Theorem
1.2(3) that $\Phi_2(C_i) \not=\emptyset$ for some $i$. If
$\Phi_2(C_i)\not= C_i$ for every $i=1, \dots ,n$, then we get an
irreducible component as in case (1). Actually, up to reordering
the curves, we can assume that $p_i \in \Phi_2(C_i)$ for $i=1,
\dots, s$. Then $G:= \langle p_1, \dots ,p_s \rangle \subseteq
\Phi_2(X)$ by Proposition 2.4. Moreover, since $\Phi_2(C_i)$ is a
finite set for every $i$, we conclude that $G$ is an irreducible
component of $\Phi_2(X)$. Now suppose that $\Phi_2(C_i)= C_i$ for
some $i$. Then up to reordering the curves we can assume that
$\Phi_2(C_i)= C_i$ for $i=1, \dots, s$ and $\Phi_2(C_i)\not= C_i$
for $i > s$. This implies that $C_i$ is a line for $i=1, \dots, s$
and a rational curve of higher degree for $i > s$. In particular
$C=\mathbb{P}^1$, i.\ e., $X$ is a rational scroll. Moreover
$G_p:=\langle p_1, \dots ,p_s \rangle \subseteq \Phi_2(X)$ by
Proposition 2.4, for every $p \in \mathbb{P}^1$. Let
$G:=\bigcup_{p\in \mathbb{P}^1} G_p$. Then $G$ is the sub-scroll
of $X$ generated by the lines $C_1, \dots C_s$. In other words,
$G$ is $\mathbb{P}(\mathcal{O}_{\mathbb{P}^1}(1)^{\oplus s}) =
\mathbb{P}^1 \times \mathbb{P}^{s-1}$ embedded in the linear span
of $C_1, \dots , C_s$ via the Segre embedding. This gives case (2)
and there are no further possibilities. The last assertions follow
from Proposition 1.4, since $G \subseteq \Phi_2(X)$.
\end{proof}
\par
Of course it may happen that a fibre $G_p$ of an irreducible
component of $\Phi_2(X)$ of type (2) is contained in a larger
component of $\Phi_2(X)$ of type (1). This happens if $C_j$ has a
flex at the point $p_j$ for some $j > s$.
\par
Now let us consider the second discriminant locus: for simplicity
we set $\mathcal{D}= \mathcal{D}_2(X,V)$ and denote by
$\mathcal{D}_G$ the component of $\mathcal{D}$ arising from an
irreducible component $G$ of $\Phi_2(X)$. Then
\[\mathcal{D}_G = \{H \in \mathbb{P}^{N \vee} \ | \ H \supseteq
\text{Osc}_x^2(X)\ \text{for any $x \in G$}\}.\] Let $G$ be an
irreducible component of $\Phi_2(X)$. We say that $G$ is of type
(1) or (2) according to the cases of Proposition 4.2. Let $G$ be
of type (1). Then, recalling that $\text{Osc}_x^2(X)$ is a fixed
$\mathbb{P}^{2n-1}$ for all $x \in G$, we conclude that the
component $\mathcal{D}_G$ is a linear $\mathbb{P}^{N-2n}$. Now
suppose that $G$ is of type (2). By Proposition 4.2,
$\text{Osc}_x^2(X)$ is a fixed linear space $T_p: =
\mathbb{P}^{2n-1}$ for $x \in G \cap f_p$. So, letting $p$ vary on
$\mathbb{P}^1$ we see that
\[\mathcal{D}_G = \bigcup_{p\in \mathbb{P}^1}\ \{H \in \mathbb{P}^{N \vee}\ |\
H \supseteq T_p \}.\] We can think of $T_p$ as
$\text{Osc}_{p_1}^2(X)$. Note that if $\Phi_2(X)\not= X$, the
tangent line to $C_n$ varies as $p$ varies on $C$. Hence for
points $x, y \in G$, lying on general distinct fibres $f_p, f_q$
of $X$ we have
\[\text{Osc}_x^2(X)=\text{Osc}_{p_1}^2(X) \not=
\text{Osc}_{q_1}^2(X)= \text{Osc}_y^2(X).\] It follows that $\dim
(\mathcal{D}_G) = N-2n+1$, $\mathcal{D}_G$ being a family of
$\mathbb{P}^{N-2n}$ parameterized by $\mathbb{P}^1$. More
precisely, we can describe the structure of $\mathcal{D}_G$ in
this way. Consider the incidence correspondence
\[\mathcal{P} =\{(p_1,H) \in C_1 \times \mathbb{P}^{N \vee}\ |\ H \supseteq
\text{Osc}_{p_1}^2(X)\}.\] Note that $\mathcal{P}$ is a
$\mathbb{P}^{N-2n}$-bundle over $C_1=\mathbb{P}^1$ via the first
projection of $C_1 \times \mathbb{P}^{N \vee}$, since
$\text{Osc}_{p_1}^2(X)$ is a $\mathbb{P}^{2n-1}$ for any $p_1 \in
C_1$. Then $\mathcal{D}_G = \pi(\mathcal{P})$, where $\pi$ is the
second projection of $C_1 \times \mathbb{P}^{N \vee}$.

\begin{example} Let $X$ be a decomposable scroll as in Section 1,
generated by lines $C_1, \dots, C_{n-1}$ and by a non-degenerate
rational curve $C_n \subset \mathbb{P}^r$, $r = r_n \geq 3$, of
degree $d$. Then $d \geq 3$ and $X \subset \mathbb{P}^N$, where
$N=2n-2+r$. According to Proposition 4.2, $G = C_1 \times
\mathbb{P}^{n-2}$, Segre embedded in $\mathbb{P}^{2n-3}=\langle
C_1, \dots C_{n-1} \rangle$. Let $p,q$ be any two distinct points
of $C=\mathbb{P}^1$. The tangent lines to $C_n$ at $p_n$ and $q_n$
generate at most a $\mathbb{P}^3$. Therefore
\[\dim \big(\langle C_1, \dots, C_{n-1},\text{Osc}_{p_n}^1(C_n),\text{Osc}_{q_n}^1(C_n)
\rangle \big) \leq 2n+1.\] Recall that the linear space above is
just the linear span $\langle \text{Osc}_{p_1}^2(X),
\text{Osc}_{q_1}^2(X)\rangle$ by Proposition 4.2. So, if $r \geq
4$, for any two distinct points $p,q \in C$ there exists a
hyperplane $H$ of $\mathbb{P}^N$ containing both
$\text{Osc}_{p_1}^2(X)$ and $\text{Osc}_{q_1}^2(X)$. Note that any
such a hyperplane corresponds to a singular point of
$\mathcal{D}_G$: actually, $\pi|_{\mathcal{P}}^{-1}(H)= \{(p_1,H),
(q_1,H)\}$. In particular, if $r \geq 4$, then
$\text{Sing}(\mathcal{D}_G)$ contains the $\mathbb{P}^2$
parameterizing the double symmetric product of $C_n$ with itself.
Now let $r=3$. If $d \geq 4$, then $C_n$ is not normal. Hence it
is the image of a rational normal curve $\overline{C} \subset
\mathbb{P}^d$ of degree $d$ via a projection from a general linear
space $T$ of dimension $d-4$. Any two tangents to $\overline{C}$
span a $\mathbb{P}^3$, so, one sees by a dimension count that in
the dual space $\mathbb{P}^{d \vee}$ there is a one dimensional
family of hyperplanes of $\mathbb{P}^d$ containing two tangent
lines to $\overline{C}$ and $T$. Projecting to $\mathbb{P}^3$ they
provide infinitely many pairs of coplanar tangent lines to $C_n$
(see also \cite[Remark 5.2]{bib5}). This being a closed condition,
implies that for any $p_n \in C_n$ there exists some other point
$q_n \in C_n$ such that the two tangent lines
$\text{Osc}_{p_n}^1(C_n)$ and $\text{Osc}_{q_n}^1(C_n)$ are
coplanar. For such a pair of points,
\[\dim \big( \langle C_1, \dots, C_{n-1}, \text{Osc}_{p_n}^1(C_n),
\text{Osc}_{q_n}^1(C_n) \rangle \big) = 2n.\] Hence $H:=
\langle\text{Osc}_{p_1}^2(X), \text{Osc}_{q_1}^2(X) \rangle$ is a
hyperplane of $\mathbb{P}^{2n+1}$ giving rise to a singular point
of $\mathcal{D}_G$. Finally, let $r=3=d$. Then $C_n$ is a twisted
cubic. Being a rational normal curve, we know that any two
distinct tangent lines to $C_n$ do not meet. Thus, for any two
distinct points $p, q \in C$ we have that
\[\langle\text{Osc}_{p_1}^2(X), \text{Osc}_{q_1}^2(X) \rangle=
\langle C_1, \dots, C_{n-1},\text{Osc}_{p_n}^1(C_n),
\text{Osc}_{q_n}^1(C_n) \rangle = \langle C_1, \dots, C_n
\rangle\] is the whole $\mathbb{P}^{2n+1}$. In fact, in this case
$\mathcal{D}_G$ is a scroll over $C$ (see Proposition 4.4 below).
\end{example}

\par
What we said in Example 4.3 when either $r \geq 4$ or $r=3$ and $d
\geq 4$ holds, ``a fortiori'', if $C_1, \dots, C_s$ are lines $(s
\geq 1$) and for some $i=s+1, \dots, n$ either $r_i \geq 4$ or
$r_i=3$ and $\deg C_i \geq 4$. Actually, also in this case there
are hyperplanes $H$ of $\mathbb{P}^N$ containing both
$\text{Osc}_{p_1}^2(X)$ and $\text{Osc}_{q_1}^2(X)$, for distinct
points $p,q$ of $C$, and any such hyperplane gives rise to a
singular point of $\mathcal{D}_G$. From now on in this Section, we
assume that
\[1 = \deg C_1 = \dots = \deg C_s < \deg C_{s+1} \leq \dots \leq \deg
C_n.\] Example 4.3 shows that $\mathcal{D}_G$ is not a scroll if
$\deg C_n \geq 4$. On the other hand, we can prove the following

\begin{proposition} Let $X$ be a decomposable scroll generated by
$C_1, \dots, C_n$, where $C_i$ is a line for $i=1, \dots, s$ and
\[2 \leq \deg C_{s+1} \leq \dots \leq \deg C_n \leq 3.\]
Let $G$ be the sub-scroll of $X$ generated by $C_1, \dots, C_s$.
Then $\mathcal{D}_G$ is a rational scroll.
\end{proposition}
\par
We need to point out some facts.

\begin{remark} Let $Y$ be the decomposable scroll generated by
$C_{s+1}, \dots , C_n$, and let $\mathbb{P}^M$ be its linear span
in $\mathbb{P}^N$. We denote by $\Sigma$ the minimal sub-scroll of
$Y$ generated by the curves $C_{s+1}, \dots , C_{n-1}$ and by $F$
any fibre of $Y$. Note that $\Sigma \cap C_n = \emptyset$, while
$\Sigma \cap F$ is a hyperplane of $F$. We have $\text{Pic}(Y)
\cong \mathbb Z^2$ and we can choose as generators the classes of
$\Sigma$ and $F$. Then:

\par (i)  any
hyperplane of $\mathbb{P}^M$ cuts $Y$ along a divisor $D$ linearly
equivalent to $\Sigma + b F$, for some integer $b > 0$. In
particular, since $\Sigma$ does not meet $C_n$ we see that
\begin{equation} \label{4.5.1}
\deg C_n = D C_n = (\Sigma + bF) C_n = b.
\end{equation}
\par (ii) For any hyperplane $H$ of $\mathbb{P}^N$ not containing $\mathbb{P}^M$
set $h:= H \cap \mathbb{P}^M$. If $h$ contains
$\text{Osc}^1_{p_i}(C_i)$ for every $i=s+1, \dots , n$ then $h$
cuts $Y$ along a divisor of the form $D = 2F_p + R$ where $F_p =
\langle p_{s+1}, \dots , p_n \rangle$ and $R$ is an effective
divisor linearly equivalent to $\Sigma + \beta F$, with $\beta
\geq 0$. Indeed, the tangent space to $Y$ at $p_i$ is
\[\text{Osc}^1_{p_i}(Y) = \langle p_{s+1}, \dots ,
\text{Osc}^1_{p_i}(C_i), \dots , p_n \rangle\] for $i=s+1, \dots ,
n$, by Remark 2.1 with $k=1$. Hence $h$ is tangent to $Y$ at all
points $p_{s+1}, \dots , p_n$. Since they are linearly
independent, this says that $h$ is tangent to $Y$ along the whole
fibre $F_p$. Thus the divisor $D$ cut out by $h$ on $Y$ is
singular at all points of $F_p$, hence the summand $2F_p$ appears
in the expression of $D$ as positive linear combination of its
irreducible components.
\end{remark}

Now we can prove Proposition 4.4.
\begin{proof} As the fibres of $\mathcal{P}$ are
mapped linearly into $\mathbb{P}^{N \vee}$ by $\pi$, it is enough
to show that the bundle projection of $\mathcal{P}$ induces a
morphism $\mathcal{D}_G \to C_1$. To do that we prove that $\pi$
is bijective, i.\ e., for any $H \in \mathcal{D}_G$, the fibre
$\pi|_{\mathcal{P}}^{-1}(H)$ consist of a single element.
Equivalently, for any pair of distinct points $p, q \in C$, there
is no hyperplane $H \subset \mathbb{P}^N$ containing both
$\text{Osc}_{p_1}^2(X)$ and $\text{Osc}_{q_1}^2(X)$. Set
\[R_p = \langle \text{Osc}_{p_{s+1}}^1(C_{s+1}), \dots , \text{Osc}_{p_n}^1(C_n)
\rangle,\quad  R_q = \langle \text{Osc}_{q_{s+1}}^1(C_{s+1}),
\dots , \text{Osc}_{q_n}^1(C_n) \rangle.\] By Proposition 4.2 we
know that
\[\text{Osc}_{p_1}^2(X)=\langle C_1, \dots, C_s, R_p \rangle, \quad
\text{Osc}_{q_1}^2(X)=\langle C_1, \dots, C_s, R_q \rangle.\] Thus
the assertion follows once we show that $\langle R_p, R_q \rangle
= \mathbb{P}^M$. By contradiction, suppose that there is a
hyperplane $h$ of $\mathbb{P}^M$ containing both $R_p$ and $R_q$.
Then, according to Remark 4.5\ (ii), $h$ cuts $Y$ along a divisor
$D = 2F_p + 2 F_q + R$, with $R$ linearly equivalent to $\Sigma +
\beta F$, for some integer $\beta \geq 0$. Then, dotting with
$C_n$ and recalling (\ref{4.5.1}), we get
\[ \deg C_n = D C_n = 4 + \beta \geq 4,\]
a contradiction.
\end{proof} \par

A further property of $\mathcal{D}_G$ is that it is degenerate in
$\mathbb{P}^{N \vee}$. In fact, for any $x \in G$,
$\text{Osc}_x^2(X)$ contains the lines $C_1, \dots, C_s$, and
hence their linear span $\Lambda:= \langle C_1, \dots, C_s
\rangle$ which is a $\mathbb{P}^{2s+1}$. By duality, this means
that $\mathcal{D}_G$ is contained in the linear subspace
$\mathbb{P}^{N-2s} \subset \mathbb{P}^{N \vee}$ parameterizing the
hyperplanes containing $\Lambda$. Moreover, $\langle \mathcal{D}_G
\rangle = \mathbb{P}^{N-2s}$.
\par
Next we want to determine the degree of $\mathcal{D}_G$ when $G$
is of type (2). To do that, recall that $\langle C_i \rangle =
\mathbb{P}^{r_i}$, and let $d_i=\deg C_i$.

\begin{proposition} Let $G$ be an irreducible component of
$\Phi_2(X)$ of type $(2)$. Then $\deg \mathcal{D}_G = 2
\sum_{i=s+1}^n (d_i-1)$.
\end{proposition}

\begin{proof} Since $G$ is of type (2) we know that $d_1= \dots =
d_s=1$ and $d_i \geq 2$ for $i > s$. Also $r_1= \dots =r_s=1$ and
$r_i \geq 2$ for $i > s$. Recalling that $N=\sum_{i=1}^n r_i +
n-1$, we note that
\[\dim \mathcal{D}_G = N-2n+1 = \sum_{i=1}^n r_i -n = \sum_{i=s+1}^n
(r_i-1).\] So $\deg \mathcal{D}_G$ is the number of elements of
$\mathcal{D}_G$ contained in a linear system $\mathcal{S} \subset
|V|$ defined by $\sum_{i=s+1}^n (r_i-1)$ linear conditions,
general enough. Choose $r_i-1$ general points in each
$\mathbb{P}^{r_i}$ for $i=s+1, \dots, n$ and call $Z_i \subset
\mathbb{P}^{r_i}$ the linear subspace they generate. Let
$\mathcal{S}$ be the linear system of hyperplanes of
$\mathbb{P}^N$ defined by the condition of passing through all
these points. Let $H \in \mathcal{D}_G$ be a hyperplane of
$\mathbb{P}^N$ not containing $\mathbb{P}^{r_i}$ for a given $i$,
$s+1 \leq i \leq n$. Then $h_i:=H \cap \mathbb{P}^{r_i}$ is a
hyperplane of $\mathbb{P}^{r_i}$ tangent to the curve $C_i$. More
precisely, if $H \supset \text{Osc}_x^2(X)$ and $x \in f_p$, then
$\text{Osc}_x^2(X) \supset \text{Osc}_{p_i}^1(C_i)$, and so $h_i$
is tangent to $C_i$ at $p_i$. On the other hand, if our $H$ is
also in $\mathcal{S}$, then, in particular, $h_i$ contains $Z_i$.
Conversely, suppose that $h_i$ is a hyperplane of
$\mathbb{P}^{r_i}$ containing $Z_i$ and tangent to $C_i$ at a
point $p_i$, and set
\[H:= \langle C_1, \dots, C_{i-1}, h_i, C_{i+1}, \dots, C_n
\rangle.\] Clearly $H$ is a hyperplane of $\mathbb{P}^N$.
Moreover, $H \in \mathcal{D}_G$, because $H$ contains all $\langle
C_j \rangle$ for $j\not=i$ and also $\text{Osc}_{p_i}^1(C_i)$.
Furthermore, $H \in \mathcal{S}$ since $H \supset Z_i$ for every
$i=s+1, \dots,n$. It thus follows that
\[\deg \mathcal{D}_G = \sum_{i=s+1}^n b_i,\]
where $b_i$ is the number of hyperplanes of $\mathbb{P}^{r_i}$
containing $Z_i$, that are tangent to $C_i$. To compute $b_i$ note
that $\dim Z_i = r_i-2$, so $Z_i$ is the axis of a pencil of
hyperplanes of $\mathbb{P}^{r_i}$. The number of hyperplanes in
this pencil that are tangent to $C_i$ is that of the ramification
points of the morphism $C_i \to \mathbb{P}^1$ defined by the
projection of $C_i$ from $Z_i$. Thus the Riemann--Hurwitz formula
tells us that $b_i=2(d_i-1)$ and this concludes the proof.
\end{proof}
\par
Relying on the above results we get the following
characterization.

\begin{theorem} Let $X \subset \mathbb{P}^N$ be a decomposable scroll
generated by $C_1, \dots, C_n$, and let $d_i=\deg C_i$, for $i=1,
\dots, n$. Suppose that $G$ is an irreducible component of type
$(2)$ of $\Phi_2(X)$. Then $\mathcal{D}_G$ is a rational normal
scroll if and only if, up to reordering the curves, $d_1= \dots
=d_s=1$ and $d_{s+1}= \dots = d_n=2$ for some $s\geq 1$.
\end{theorem}

\begin{proof} As $G$ is of type (2), we can assume that $d_1= \dots
=d_s=1$ for some $s \geq 1$ and $d_i \geq 2$ for $i \geq s+1$, by
Proposition 4.2. As we noted, $\mathcal{D}_G$ has dimension
$N-2n+1$ and is non-degenerate in $\mathbb{P}^{N-2s}$. Thus,
recalling Proposition 4.6, the inequality $\deg \mathcal{D}_G \geq
\text{codim}\mathcal{D}_G + 1$ becomes
\begin{equation}\label{4.7.1} 2 \sum_{i=s+1}^n (d_i-1) \geq 2(n-s). \end{equation}
Note that this is an equality if and only if
\begin{equation}\label{4.7.2} d_i=2 \quad \text{for} \quad i=s+1 \dots, n. \end{equation}
So, if $\mathcal{D}_G$ is a rational normal scroll, then
(\ref{4.7.2}) holds. On the other hand, if (\ref{4.7.2}) holds
then we know that $X$ is a rational scroll, by Proposition 4.4,
and then equality in (\ref{4.7.1}) says that it is normal.
\end{proof}

\section{A general lower bound}

In \cite[Theorem A]{bib6} it is shown that the highest
inflectional locus of an indecomposable linearly normal elliptic
scroll of invariant $e=-1$ is empty. By adapting the argument used
in \cite{bib6} we can locate the highest inflectional locus of an
elliptic indecomposable scroll of invariant $e=0$, which is
linearly normally embedded. Let $C$ be a smooth curve of genus $1$
and let $S=\mathbb{P}(\mathcal{E})$, where $\mathcal{E}$ is the
holomorphic rank-2 vector bundle on $C$ defined by the non-split
extension
\begin{equation}\label{5.0.1} 0 \to \mathcal{O}_C \to \mathcal{E} \to \mathcal{O}_C
\to 0. \end{equation} Let $\pi:S \to C$ be the ruling projection
and denote by $C_0$ the tautological section on $S$. Let $\delta
\in \text{Div}(C)$ be a divisor of degree $\deg \delta = m+1 \geq
3$ and set $L: = \mathcal{O}_S(C_0+\pi^* \delta)$. Note that $L$
is very ample, because $\deg \delta \geq e+3$ \cite[Ex.\ 2.12(b),
p.\ 385]{bib4} and the morphism given by $|L|$ embeds $S$ as a
linearly normal scroll of degree $2m+2$ in $\mathbb{P}^N$, where
$N=2m+1$ (note that $m=\frac{N-1}{2}=\lbrack \frac{N-1}{2}
\rbrack$). Let $x \in S$. By \cite[(1.$0_m$)]{bib6} we have
\begin{equation}\label{5.0.2} \dim \big( \text{Osc}_x^m(S) \big) = N-1-
\dim(|L-(m+1)x|) = 2m - \dim(|L-(m+1)x|). \end{equation}
On the
other hand, by \cite[Remark 1.2]{bib6} we know that
\begin{equation}\label{5.0.3} |L-(m+1)x| = m f_x + |L - mf_x - x|,
\end{equation} where $f_x$ is the fibre through $x$. Note that the line
bundle $L\otimes \mathcal{O}_S(-mf_x) =
\mathcal{O}_S(C_0+\pi^*(\delta - m \pi(x))$ is not necessarily
spanned, because $\deg(\delta - m\pi(x))=m+1-m=1 < e+2$ \cite[Ex.\
2.12(a), p.\ 385]{bib4}. We have that
\begin{equation} \label{5.0.4}
\dim(|L-mf_x-x|) =\dim(|L-mf_x|)-1\ \ \text{if and only if
$L\otimes \mathcal{O}_S(-mf_x)$}\ \text{is spanned at $x$}.
\end{equation}
Now, twisting (\ref{5.0.1}) by $\mathcal{O}_C(\delta - m \pi(x))$
we immediately see that
\[h^0(L-mf_x) = h^0(\mathcal{E}(\delta - m\pi(x)) = 2
h^0(\mathcal{O}_C(\delta - m\pi(x)) = 2.\] Hence (\ref{5.0.4})
gives
\[\dim(|L-mf_x-x|)=0 \qquad \text{if and only if $L\otimes \mathcal{O}_S(-mf_x)$ is spanned at
$x$}\] and taking into account (\ref{5.0.3}) and (\ref{5.0.2}) we
get
\[\dim\big( \text{Osc}_x^m(S) \big) = 2m \qquad \text{if and only if $L\otimes \mathcal{O}_S(-mf_x)$ is
spanned at $x$.}\] This proves the following

\begin{proposition} Let $S \subset \mathbb{P}^{2m+1}$ be a linearly
normal surface scroll over an elliptic curve $C$, defined by an
indecomposable vector bundle as in $(\ref{5.0.1})$, and let $L$ be
the hyperplane bundle. Then $x \in \Phi_m(S)$ if and only if the
line bundle $L\otimes \mathcal{O}_S(-mf_x)$ is not spanned at $x$,
where $f_x$ is the fibre of $S$ through $x$.
\end{proposition}

\par
Now let $S \subset \mathbb{P}^N$ be any surface scroll. Though $S$
can be not decomposable, according to Example 3.2 it seems natural
to ask whether, under some assumption, we can get a global lower
bound for the dimension of $\text{Osc}^k_x(S)$, i.\ e., a lower
bound holding at every point $x \in S$, bigger than $3$. We
determine such a lower bound, depending on $k$, under certain
assumptions on the linear system (not necessarily complete) giving
rise to the embedding.

In the following we use the same notation as in \cite{bib6}.

\begin{theorem} Let $S \subset \mathbb{P}^N$ be a surface scroll
embedded by $|V|$, where $V \subseteq H^0(S,L)$,
$L=\mathcal{O}_{\mathbb{P}^N}(1)|_S$, and suppose that $N \geq
2k$. Let $x \in S$ and denote by $f_x$ the fibre of $S$ through
$x$. If $|V-tf_x|$ is very ample for every non-negative integer $t
\leq k-2$, then
\[\dim\big(\text{\rm{Osc}}_x^k(S)\big) \geq k+2 \qquad \text{for}
\ k \geq 3.\]
\end{theorem}

Note that this global lower bound is the same holding for rational
normal scrolls, as shown by Example 3.1.

\begin{proof} First let us prove, by induction, that
\begin{equation} \label{5.2.1}
\dim\big(\text{Osc}_x^k(S)\big) \geq k+1 \qquad \text{\rm{for any
$k \geq 2$}}. \end{equation} For $k=2$ this comes from
\cite[Theorem B]{bib6} (noting that $N \geq 4$ is enough in the
proof). So let $k \geq 3$ and set $\mathcal{L}=L-(k-2)f_x$ and
$|W|=|V-(k-2)f_x|$. Note that $|W|$ is very ample by assumption
and that $S$ embedded by $|W|$ is also a scroll. Actually, for
every fibre $f$ of $S$ we have
\[\mathcal L f = (L-(k-2)f_x) f = L f = 1.\] Thus \cite[Lemma (1.4)
and Lemma (1.5)]{bib6} imply the following facts. The linear
system $|W - f_x|$ is base-point free and if $\varphi_x:S \to
\mathbb{P}^M$ denotes the associated morphism then one of the
following conditions holds:
\begin{enumerate}
\item[{\rm i)}] every fibre of $\varphi_x$ intersects any fibre
$f$ of $S$ at a finite set, \item[{\rm ii)}] $(S,\mathcal{L}) =
(\mathbb{P}^1 \times \mathbb{P}^1, \mathcal{O}_{\mathbb{P}^1
\times \mathbb{P}^1}(1,1))$ and $W = H^0(S,\mathcal{L})$.
\end{enumerate}
If ii) holds, then $L=\mathcal{O}(1,k-1)$, hence $h^0(L)=2k$,
which implies that \[N = \dim(|V|) \leq \dim(|L|) = 2k-1,\] but
this contradicts our assumption that $N \geq 2k$. Therefore
condition i) holds. Now suppose, by contradiction, that
(\ref{5.2.1}) is not true, i.\ e., $\dim\big( \text{Osc}_x^k(S)
\big) \leq k$. Because $|V -(k-3)f_x|$ is very ample, by induction
we know that $\dim \big( \text{Osc}_x^{k-1}(S) \big) \geq k$. So,
due to the obvious inclusion $\text{Osc}_x^k(S) \supseteq
\text{Osc}_x^{k-1}(S)$ we conclude that
\[\text{Osc}_x^k(S) = \text{Osc}_x^{k-1}(S).\] Equivalently, this
says that \[|V-(k+1)x| = |V-kx|.\] This in turn, according to
\cite[Remark (1.2)]{bib6}, implies the equality
\[\dim(|V-kf_x-x|) = \dim(|V-(k-1)f_x-x|).\] Hence \[f_x \subseteq
\text{Bs}(|V-(k-1)f_x-x|) = \text{Bs}(|W-f_x-x|) =
\varphi_x^{-1}(\varphi_x(x)).\] But this contradicts condition i).
To conclude the proof we show that equality cannot occur in
(\ref{5.2.1}) for $k\geq 3$. First of all, since $|V-tf_x|$ is
very ample for all $t \leq k-2$, by applying \cite[Remark
(1.7)]{bib6} inductively we see that
\begin{equation} \label{5.2.2}
\dim(|W|) = \dim(|V|) - 2(k-2) = N - 2(k-2). \end{equation} Note
that
\[\dim(|V-(k+1)x|)= \dim(|W-3x|)\] by \cite[Remark (1.2)]{bib6}.
Due to (\ref{5.2.1}) and the assumption $N \geq 2k$, $S$ is
embedded by $|W|$ as a scroll in a projective space of dimension
$\geq 4$; hence $|W-3x| \not= |W-2x|$ by \cite[Theorem B]{bib6}.
Since $|W|$ is very ample this says that $\dim(|W-3x|)<
\dim(|W|)-3$. Thus, recalling the equality
\[\dim(|V-(k+1)x|)+\dim \big( \text{Osc}_x^k(S) \big)=N-1,\] we get
\begin{equation} \label{5.2.3}
\dim \big( \text{Osc}_x^k(S) \big) > N+2-\dim(|W|).
\end{equation}
Finally, combining (\ref{5.2.2}) with (\ref{5.2.3}) and assuming
equality in (\ref{5.2.1}) gives $k \leq 2$. This completes the
proof.
\end{proof}

\begin{acknowledgement} During the preparation of this paper the first author has been
supported by the MUR of the Italian Government in the framework of
the PRIN ``Geometry on Algebraic Varieties'', and the second
author by the projects BFM2003-03917/MATE (Spanish Ministry of
Education) and Santander/UCM PR27/05-13876. The first author would
also like to thank the GNSAGA-INDAM and Azione Integrata
Italia--Spagna IT200 for support received at an early stage of
this research. Both authors are grateful to the University of
Milan for financial support.
\end{acknowledgement}

\end{document}